\input amstex

\documentstyle{amsppt}
\magnification=\magstep1

\pageheight{7in}
\TagsOnRight
\NoRunningHeads
\NoBlackBoxes

\def\tt#1{\medskip\leftline{\bf #1}\medskip\flushpar\hglue0pt}
\def\tit#1{\par\medskip\noindent{\bf #1}\par\smallskip\noindent\hglue0pt}

\def\ci#1{_{{}_{\ssize #1}}}
\def\cci#1{_{{}_{\sssize #1}}}
\def\ui#1{^{{}^{{}_ {\,#1} }}}

\def\R{\Bbb R}
\def\C{\Bbb C}

\def\LHS{\operatorname{LHS}}
\def\RHS{\operatorname{RHS}}
\def\BMO{\operatorname{BMO}}
\def\dist{\operatorname{dist}}
\def\dim{\operatorname{dim}}
\def\Vol{\operatorname{Vol}}
\def\vol{\operatorname{Vol}}

\def\x{\bold x}
\def\y{\bold y}

\redefine\ge{\geqslant}
\redefine\le{\leqslant}

\def\dd{\partial}

\def\la{\lambda}
\def\La{\Lambda}
\def\f{\varphi}
\def\e{\varepsilon}
\def\wt{\widetilde}
\def\wh{\widehat}

\def\a{\alpha}
\def\d{\delta}

\def\<{\langle}
\def\>{\rangle}

\topmatter

\title
The geometric Kannan-Lov\'asz-Simonovits lemma, dimension-free estimates
for volumes of sublevel sets of polynomials, and distribution of
zeroes of random analytic functions.
\endtitle

\thanks This research was partially supported by the United States
- Israel Binational Science Foundation \endthanks

\author
F. Nazarov, M. Sodin, A. Volberg
\endauthor

\address (F.N.:) Department of Mathematics, Michigan State
University, East Lansing, MI 48824, U.S.A. 
\endaddress
\email fedja\@math.msu.edu \endemail

\address (M.S.:) School of Mathematical Sciences, Tel Aviv University,
Ramat Aviv, 69978, Israel
\endaddress
\email sodin\@post.tau.ac.il \endemail

\address (A.V.:) Department of Mathematics, Michigan State
University, East Lansing, MI 48824, U.S.A. 
\endaddress
\email volberg\@math.msu.edu \endemail

\abstract
The goal of this paper is to attract the attention of the reader to a
simple dimension-free geometric inequality that can be proved using the
classical needle decomposition technique. This inequality allows us to
derive
sharp dimension-free estimates for the distribution of values
of polynomials in convex subsets in $\R^n$ in a simple and elegant way.
Such estimates, in their turn, lead to a surprising result about the
distribution of zeroes of random analytic functions; informally speaking,
we show that for simple families of analytic functions, there exists a
``typical" distribution of zeroes such that the
portion of the family occupied by the functions
whose distribution of zeroes deviates from
that typical one by some
fixed amount is about
$\operatorname{Const}\exp\{-\operatorname{size\ of\ the\ deviation}\}$.

The paper is essentially self-contained.
When choosing the style, we tried to make it an enjoyable reading for both a
senior undergraduate student and an expert.

As to the standard question ``What is new in the paper?" one is
supposed to address in the abstract, we believe that the answer to it
is a function of {\it two} variables, the first being ``what is written" and
the second being ``who is reading". Since we have no knowledge of the value 
of
the second variable, we can only give the range of answers with the first
variable fixed. We believe that for the targeted audience it will be
the standard range $[\operatorname{Nothing},
\operatorname{Everything}]$ (with both endpoints included).

\endabstract

\endtopmatter

\head
$\S 1$. The Geometric Kannan-Lov\'asz-Simonovits Lemma
\endhead

\par\noindent By this name we will call the following

\tit{Proposition:}
{\it Let $\Cal F$ be a compact convex subset of $\R^n$
with non-empty interior, and let $\Cal E\subset \Cal F$ be a closed set.
Let $\la >1$, and let 
$$
\Cal E\ci{\la,\Cal F}:= \Bigl\{
\x\in \Cal E:\, \text{for every interval $\Cal J$ such that }\x\in
\Cal J\subset \Cal F,\
\frac{|\Cal E\cap \Cal J|}{|\Cal J|}\ge \frac{\la-1}\la
\Bigr\}.
$$
Then}
$$
\frac{\vol(\Cal E\ci{\la,\Cal F})}{\vol(\Cal F)}\le
\left[
\frac{\vol(\Cal E)}{\vol(\Cal F)}
\right]^\la.
$$
\tit{Remark:}
In the definition of the ``core'' $\Cal E\ci{\la,\Cal F}$ it
is enough to consider only
the intervals $\Cal J$ that have $\x$ as one of their endpoints. Indeed,
if $\x$ is an interior point of $\Cal J$, and the condition
$\dfrac{|\Cal E\cap \cdot|}{|\cdot|}\ge \dfrac{\la-1}\la$ is satisfied for
each
of the two subintervals into which $\x$ splits $\Cal J$,
then it is satisfied for the entire interval $\Cal J$.

\tit{Proof of the Geometric KLS Lemma:}
Consider first the following special case: let
$\Cal L\subset \R^n$ be a line and
let $\Bbb P$ be the orthogonal projection to $\Cal L$. Let
$I=\Bbb P \Cal F\subset \Cal L$.
Assume that $\Cal E=\{\x\in \Cal F\,:\,\Bbb P\x\in E\}$ where $E$ is some
closed subset of $I$.

\tit{Claim:}
$$
\Cal E\ci{\la,\Cal F}=\Cal F\cap \Bbb P^{-1} E\ci{\la,I}.
$$
Informally, this means that the set $\Cal E\ci{\la,\Cal F}$
is determined by its projection onto the line $\Cal L$
(as the maximal subset of $\Cal F$ with given projection)
and that this projection is
$$
E\ci{\la,I} := \Bigl\{
x\in E:\, \text{for every interval $J$ such that }x\in J
\subset I,\
\frac{|E\cap J|}{|J|}\ge \frac{\la-1}\la
\Bigr\}.
$$
Speaking even more vaguely, one may say that the statement
of the geometric KLS lemma
for such ``simple'' sets is ``essentially one-dimensional".

\tit{Proof of the claim:}
Since this claim is a simple exercise in geometry, we shall
present only the proof of that part of it that we really need,
namely, of that the set on the
left is contained in the set on the right.
Suppose that $\x\in \Cal E$ and $x=\Bbb P \x\notin E\ci{\la,I}$.
Then there exists an interval $J\subset I$ such that $x$ is one of
its endpoints and $\dfrac{|J\cap E|}{|J|}<\dfrac{\la-1}{\la}$.
Let $y\in I$ be the other end of $J$. There exists a point $\y\in \Cal F$
such that $y=\Bbb P \y$. Since $\Cal F$ is convex,
the entire interval $\Cal J=\x\y$
is contained in $\Cal F$. It is easy to check that
$\dfrac{|\Cal J\cap \Cal E|}{|\Cal J|}=
\dfrac{|J\cap E|}{|J|}<\dfrac{\la-1}{\la}$ and thereby
$\x\notin \Cal E\ci{\la,\Cal F}$. $\square$

\medskip
Once the claim has been proved, we are ready to reformulate the
statement of the geometric KLS lemma for this special case as
a one-dimensional problem. Let $f(x)$ ($x\in I$)
be the $(n-1)$-dimensional volume of the
cross-section of the convex set $\Cal F$ by the hyperplane orthogonal to the
line $\Cal L$ and containing the point $x$.
We have
$$
\Vol(\Cal E\ci{\la,\Cal F})=
\int_{E\ci{\la,I}}f(x)\,dx;
$$
$$
\Vol(\Cal F)=\int_{I}f(x)\,dx;
$$
and
$$
\Vol(\Cal E)=\int_{E} f(x)\,dx.
$$
Using these three formulae,
we see that the statement of the geometric KLS lemma for
our special case can be rewritten as
$$
\frac{\int_{E\cci{\la,I}} f}{\int_{I}f}\le
\left[
\frac{\int_{E} f}{\int_{I}f}
\right]^\la.
$$
The best thing one can hope for is that this inequality is valid
for an {\it arbitrary} non-negative continuous
function $f$ and an arbitrary set $E\subset I$. It doesn't take a long
time to see that it is not the case, so
the next natural question to ask is ``What is so special about the
functions that express the volumes of cross-sections of convex bodies?".
The answer is given by the classical
Brunn-Minkowski theorem, one of several equivalent formulations of which
is that the function $f(x)^{\frac{1}{n-1}}$ is concave, i.e.,
$$
f\bigl(tx+(1-t)y\bigr)^{\frac{1}{n-1}}\ge
t f\left(x\right)^{\frac{1}{n-1}}+
(1-t) f\left(y\right)^{\frac{1}{n-1}}\qquad\text{ for all }x,y\in I,\
t\in[0,1].
$$
This property is for each $n$ stronger than and for large $n$ almost
equivalent to {\it logarithmic concavity} of the function $f$,
i.e., to the inequality $f\bigl(tx+(1-t)y\bigr)\ge
f(x)^{t}\,f(y)^{1-t}$.
Thus, our special case is covered by the following

\tit{Lemma:}
{\it Let $I\subset \R$ be an interval and let $f:I\to[0,+\infty)$
be a logarithmically concave function that does not vanish at interior
points
of $I$.
Let $E\subset I$ be a measurable set. Fix $\la>1$ and define
$$
E\ci{\la,I} := \Bigl\{
x\in E:\, \text{for every interval $J$ such that }x\in J\subset I,\
\frac{|E\cap J|}{|J|}\ge \frac{\la-1}\la
\Bigr\}.
$$
Then}
$$
\frac{\int_{E\cci{\la,I}} f}{\int_{I}f}\le
\left[
\frac{\int_{E} f}{\int_{I}f}
\right]^\la \,.
$$

\smallskip
If at this stage the reader has the feeling that, once formulated,
this statement requires only some routine techniques he already knows
to prove it, he is probably right. We offer such a reader to
try to prove the lemma by himself before reading our proof in the
Appendix in the hope that he might be able to come up with a nicer
proof than that of ours, which, though completely natural, lacks in
elegance.

\medskip
Our next task will be to reduce the full statement of the
Geometric KLS lemma to this special case. We will do it
using the {\it classical needle
decomposition}.

First of all,
let us remind/tell the reader what the classical needle decomposition is.
Given a compact convex body $\Cal F\in\R^n$ and $\d>0$,
we can perform the following
construction. Take any $2$-dimensional plane
$\Cal K$ that intersects $\Cal F$ and choose
a $\d$-net in the set  $\Cal F\cap \Cal K$. For each point in this
$\d$-net, take an
($n-2$)-dimensional plane that is orthogonal to $\Cal K$ and intersects
$\Cal K$ at the
corresponding point. Clearly, for any $2$-dimensional plane $\Cal K'$
sufficiently close to $\Cal K$
in some natural metric
\footnote
{
One possible way to introduce a ``natural distance" between two planes
$\Cal K_1$ and $\Cal K_2$ of
the same dimension is the following. Consider all isometric motions of
$\R^n$
that map $\Cal K_1$ to $\Cal K_2$. Every such isometric motion is of the
kind
$x\to Ux+a$ where $a\in \R^n$ and $U$ is a unitary operator.
Define $\dist(\Cal K_1,\Cal K_2):=\inf(\|U-I\|+|a|)$.
To check the axioms of distance is left to the reader as an exercise.
}
these planes are transversal to $\Cal K'$ and their
intersections with $\Cal K'$ form a $2\d$-net in $\Cal K'\cap \Cal F$.
Therefore, since the set of all
$2$-dimensional planes intersecting $\Cal F$ is compact in any natural
metric,
we can find finitely many ($n-2$)-dimensional spaces
$\Cal M_1,\dots,\Cal M\ci N$ such
that for every $2$-dimensional plane $\Cal K$ intersecting
$\Cal F$, the set of points
at which $\Cal K$ is intersected by those of the planes
$\Cal M_1,\dots, \Cal M\ci N$ that
are transversal to it, forms a $2\d$-net in $\Cal K\cap \Cal F$.

Carry out the following algorithm:
\tit{Step 1:}
Choose a hyperplane $\Cal H\supset \Cal M_1$.
It splits $\Cal F$ into two compact convex
subsets $\Cal F\ui{ +}$  and $\Cal F\ui{ -}$.

\tit{Step 2:}
Choose a hyperplane $\Cal H\ui{ +}\supset \Cal M_2$ and split
$\Cal F\ui{ +}$ into
$2$ compact convex
subsets (one of which may be empty) by this hyperplane.
Then choose a hyperplane $\Cal H\ui{ -}\supset \Cal M_2$ and split
$\Cal F\ui{ -}$ into
$2$ compact convex subsets by that hyperplane.

\centerline {\bf\vdots}

\tit{Step $k$:}
After completion of Step $k-1$, we have a decomposition of $\Cal F$ into
$2^{k-1}$
subsets. Split each of those subsets into two smaller ones by a hyperplane
containing $\Cal M_k$ (so, Step $k$ consists of $2^{k-1}$ substeps).
\bigskip
After completing all $N$ steps in this algorithm, we obtain a decomposition
of $\Cal F$ into $2^N$ compact convex subsets
$\Cal F_j$ some of which may be empty.
\tit{Definition:}
Let $\gamma>0$. A convex set $\Cal F$ is called a $\gamma-needle$ if there
is a
line in $\R^n$ such that the distance from every point of
$\Cal F$ to this line is not greater than $\gamma$.
\medskip
\tit{Claim:}
{\it Every set $\Cal F_j$ is an $8\d$-needle.}
\tit{Proof of the claim:}
Let us first show that for every two-dimensional plane $\Cal K$,
the set $\Cal F_j\cap K$ contains no disk $\Cal D$ of radius $2\d$.
Indeed, otherwise there
would exist an ($n-2$)-dimensional plane $\Cal M_k$ transversal
to $\Cal K$ such that
$\Cal M_k$ intersects $\Cal K$ at some point inside the disk $\Cal D$. But
then
the set $\Cal F_j$ cannot be contained entirely in any half-space
bounded by any hyperplane containing $\Cal M_k$. On the other hand,
Step $k$ provides such a half-space and we get a contradiction.

Now, let $\bold a$ and $\bold b$
be the endpoints of the longest interval contained in
$\Cal F_j$. Note that for every point
$\bold c\in \Cal F_j$, the angles $\wh {\bold a}$ and $\wh{\bold b}$
of the triangle $\bold a\bold b\bold c$ are less than $\frac{\pi}{2}$.

If $\dist(\bold a, \bold b)\le 8\d$, then $\Cal F_j$ lies in a
$8\d$-neighborhood of
any line containing the point $\bold a$.
Otherwise, consider any point $\bold c\in \Cal F_j$.
If the distance from the point $\bold c$ to the line
$\bold a\bold b$ is greater than $8\d$,
then the triangle $\bold a\bold b\bold c$
contains a rectangle both sides of which are greater
than $4\d$ and, thereby, a disk of radius $2\d$, which is impossible. Thus,
$\Cal F_j$ lies entirely in a $8\d$-neighborhood of the line $\bold a\bold
b$. $\square$

\medskip
This construction can be used (and/or generalized)
in many different ways. Since we are not after that extremely elusive
thing known by the name ``full generality" in this note, we shall
only show how this construction can be used to fit our purposes.
For other usages see the papers [ND1] by Gromov-Milman,
[ND2] by Lov\'asz-Simonovits, and [ND3] by
Kannan-Lov\'asz-Simonovits where this elementary idea
was developed into a
powerful tool in ``high-dimensional" geometry, especially in the
study of isoperimetric inequalities.

The only freedom we have in the algorithm described above is the choice of
the hyperplanes containing given ($n-2$)-dimensional planes. That is one
degree
of freedom at each substep and we can use it to ``solve one equation".

Now take some small $\d>0$. Let $\wt{\Cal E}:=
\{x\in \Cal F\,:\,\dist(x,\Cal E)\le 16\d\}$
and let $\a=\dfrac{\Vol(\wt{\Cal E})}{\Vol(\Cal F)}$.
Let us look at the first step
in the needle decomposition construction.
To choose a hyperplane $\Cal H\supset \Cal M_1$ is the
same as to choose a unit vector
$\bold v\bot \Cal M_1$ (the unit vector orthogonal to $\Cal H$).
Since $\dim \Cal M_1=n-2$, the
set of such vectors $\bold v$ is a unit circumference. Let's adopt
the natural agreement that
$\Cal F\ui{ +}=\Cal F\ui{ +}(v)$ is the part of $\Cal F$
contained in the half-space that lies in the
direction of the vector $\bold v$ from $\Cal H$, i.e.,
$$
\Cal F\ui{ +}(\bold v)=
\{\x\in \Cal F\,:\,\langle \x-\y,\bold v\rangle\ge 0\text{ for all }\y
\in \Cal H\},
$$
and
that $\Cal F\ui{ -}$ is the other part.

Suppose that $\Vol(\wt{\Cal E}\ui{ +}(\bold v))>\a \Vol (\Cal F\ui{ +}(\bold
v))$.
Then, obviously,
$$
\Vol(\wt{\Cal E}\ui{ +}(-\bold v))=\Vol(\wt{\Cal E}\ui{ -}(\bold v))
<\a\Vol(\Cal F\ui{ -}(\bold v))=\a\Vol(\Cal F\ui{ +}(-\bold v)).
$$
Thus, since the continuous function
$\bold v\mapsto \Vol(\wt{\Cal E}\ui{ +}(\bold v))-\a\Vol(\Cal F\ui{
+}(\bold
v))$ attains
both positive and negative values on the unit circumference, it must vanish
somewhere, i.e., there exists a hyperplane $\Cal H\supset \Cal M_1$ such
that
$\Vol(\wt{\Cal E}\ui{ +})=\a\Vol(\Cal F\ui{ +})$
(this is exactly that ``one equation" we
solve using one degree of freedom we have in Step~1).
Obviously, for such a hyperplane, we also have
$\Vol(\wt{\Cal E}\ui{ -})=\a\Vol(\Cal F\ui{ -})$.
It is easy to check that two other possible assumptions
$\Vol(\wt{\Cal E}\ui{ +}(\bold v))<\a \Vol (\Cal F\ui{ +}(\bold v))$
and $\Vol(\wt{\Cal E}\ui{ +}(\bold v))=\a \Vol (\Cal F\ui{ +}(\bold v))$
result in the same conclusion.

Making an analogous choice during each (sub)step, we shall arrive at the
decomposition of $\Cal F$ into $8\d$-needles $\Cal F_j$ such that the
volumes of
the corresponding parts $\wt{\Cal E}_j=\wt{\Cal E}\cap \Cal F_j$ of
the set $\wt{\Cal E}$ satisfy $\Vol(\wt{\Cal E}_j)=\a\Vol(\Cal F_j)$.

Let $\Cal L_j\subset \R^n$ be some line in whose
$8\d$-neighborhood the set $\Cal F_j$
is contained. Let $\Bbb P_j$ be the orthogonal projection onto $\Cal L_j$.
Let $I_j=\Bbb P_j \Cal F_j$. At last, let $\Cal E_j=\Cal E\cap F_j$.
Denote by $\Cal G_j$ the maximal subset of
$\Cal F_j$ whose orthogonal projection to
the line $\Cal L_j$ coincides with that of $\Cal E_j$.
Formally, it means that $\Cal G_j=\Cal F_j\cap \Bbb P_j^{-1}(\Bbb P_j
\Cal E_j)$. Clearly, $\Cal E_j\subset \Cal G_j\subset \wt{\Cal E}_j$. We
have
$$
\Cal F_j\cap \Cal E\ci{\la,\Cal F}
\subset (\Cal E_j)\ci{\la,\Cal F_j}
\subset (\Cal G_j)\ci{\la,\Cal F_j}
$$
Applying the special case of the geometric KLS lemma to the sets
$\Cal G_j$ and $\Cal F_j$ and recalling that
$\Vol(\Cal G_j)\le \Vol(\wt{\Cal E}_j)= \a\Vol(\Cal F_j)$, we obtain
$$
\Vol(\Cal F_j\cap \Cal E\ci{\la,\Cal F})
\le \a^{\la}\Vol(\Cal F_j).
$$
Adding these estimates for all $j$, we arrive at the inequality
$\Vol(\Cal E\ci{\la,\Cal F})\le \a^{\la}\Vol(\Cal F)$ or, equivalently,
$$
\frac{\Vol(\Cal E\cci{\la,\Cal F})}{\Vol(\Cal F)}\le
\left\{\frac{\Vol(\wt{\Cal E})}{\Vol(\Cal F)}\right\}^{\la}.
$$
Now, to finish the proof,
it remains only to note that $\Vol(\wt{\Cal E})\to\Vol(\Cal E)$ as $\d\to
0$. $\square$

\medskip
If the reader wants to understand this proof better and to
see how neatly the needle decomposition works, we recommend
him to consider the convex set
$\Cal F=\{\x=(x_1,x_2)\in\R^2\,:\,|\x|\le 1\}\subset\R^2$ with subsets
$\Cal E\cci 1=\{\x\in \Cal F\,:\, x_2\ge 0\}$ and
$\Cal E\cci 2=\{\x\in \Cal F\,:\, |\x|\ge r\}$ ($0<r<1$),
and draw all the corresponding pictures and write
the corresponding inequalities for these two cases.

\tt{Remark:}
An expert may observe here that, instead of volume,
we might consider an arbitrary logarithmically concave measure $\mu$
in $\R^n$, i.e., a measure of the kind $d\mu(\x)=p(\x)\,d\x$ where
the density $p:\R^n\to [0,+\infty)$ is a logarithmically concave function.
(as above, we call $p$ logarithmically concave if
$p(t\x+(1-t)\y)\ge p(\x)^t p(\y)^{1-t}$ for all $\x,\y\in\R^n$,
$0\le t\le 1$\,). When $p\equiv 1$, we get the volume. Another
interesting example coming from the probability theory is
$p(\x)=(2\pi)^{-\frac{n}{2}}e^{-\frac{|\x|^2}{2}}$, the density of
the {\it standard Gaussian distribution} in $\R^n$.

A version of the Brunn-Minkowski theorem asserts that the class of
logarithmically concave measures is closed under projections of $\R^n$ to
affine subspaces (see, e.g., [ND1]). This allows us to extend the
inequality of the geometric KLS lemma verbatim to
arbitrary finite logarithmically concave measures:
$$
\frac{\mu(\Cal E\ci{\la,\Cal F})}{\mu(\Cal F)}\le
\left[
\frac{\mu(\Cal E)}{\mu(\Cal F)}
\right]^\la
$$
for every convex set $\Cal F$ with $0<\mu(\Cal F)<+\infty$,
closed subset $\Cal E\subset\Cal F$ and $\la>1$.

On the other hand, our whole point was to re-emphasize the geometric
nature of the Lovasz-Simonovits localization technique and to somewhat
counterbalance the tendency to present the needle
decomposition as a statement about two (or four) integrals rather than a
geometric partition algorithm. So, we preferred to use ``purely geometric"
terminology and to restrict ourselves to ``volumes" and ``convex sets"
in the main text. At last, it may be worth mentioning that the class of 
logarithmically
concave measures is only marginally wider than the class of convex sets:
every logarithmically concave measure in $\R^n$ can be obtain as the limit
when $m\to\infty$ of projections to $\R^n$ of volumes of convex sets in
$\R^m\supset \R^n$. This allows one to extend many statements about
convex sets to the case of logarithmically concave measures more or less
automatically. The Gaussian measure, for example, can be viewed as the
limit of projections of the volume measures of balls in large dimensions.

\head
$\S 2$. Dimension-free estimates for volumes of sublevel sets polynomials
\endhead

Let us start with recalling the classical $1$-dimensional
\tit{Remez inequality:}
{\it Let $P$ be a polynomial of degree $d$ in $\R^1$.
Then for every interval $J\subset \R^1$ and for every measurable subset
$E\subset J$,}
$$
\max_J |P|\le \left[
\frac{A|J|}{|E|}
\right]^d\sup_E |P|
\tag R
$$
{\it where $A>0$ is an absolute constant (whose best possible value
is $A=4$).}

\smallskip
The proof (with a worse constant $A=2e$) follows by a
straightforward application of the Lagrange interpolation formula with
$d+1$ nodes on $E$ spaced by at least $\frac{|E|}{d}$.
The sharp constant can obtained by a Markov-type ``moving zeroes" argument,
which shows that the worst case is attained when $E$ is a sub-interval of
$J$ with a common end-point with $J$ and $P$ is the (properly renormalized)
Chebyshev polynomial.

There is no hope for a dimension-free Remez type inequality with the
$L^\infty$ norm on the left hand side. This can be
already seen when $\Cal F$ is a unit ball in $\R^n$ and
$P(\x)=1-|\x|^2$ (the reason is that, for large $n$, most of the volume
of $\Cal F$ is concentrated in a small neighborhood of the unit sphere
where $P(\x)$ is very small).
So, we have to confine ourselves to weaker distribution estimates.

Observing that a restriction of a polynomial of degree $d$ to any line in
$\R^n$ is again a polynomial of degree (not exceeding) $d$ and
combining the one-dimensional Remez inequality with
the geometric KLS lemma, we obtain the following
\tit{Comparison lemma:}
{\it Let $P$ be a polynomial of degree $d$ in $\R^n$, and let $\Cal F$
be a convex compact set of volume one. Then for any } $c>0$, $\la\ge
1$
$$
\vol\{\x\in \Cal F\,:\, |P(\x)|\ge (A\la)^d c \}\le
\bigl[\vol\{\x\in \Cal F\,:\,|P(\x)|\ge c\}\bigr]^{\la}.
$$
\tt{Proof:}
If $P$ is constant, the estimate is trivial. Otherwise, let
$\Cal E=\{\x\in \Cal F\,:\,|P(\x)|\ge c\}$. For each
$\x\notin \Cal E\ci{\la,\Cal  F}$, we can find an interval
$\Cal J\subset \Cal F$ containing the point $\x$ and such that
the length of the set $\Cal J\setminus \Cal E$ is at least
$\la^{-1}|\Cal J|$. Then, according to the Remez inequality,
$$
|P(\x)|\le \max_\Cal J |P|\le \left[
\frac{A|\Cal J|}{|\Cal J \setminus \Cal E|}
\right]^d\sup_{\Cal J \setminus \Cal E} |P|\le (A\la)^d c
$$
and, thereby,
$\{\x\in \Cal F\,:\, |P(\x)|> (A\la)^d c \}\subset \Cal E\ci{\la,\Cal F}$.
It remains to observe that the strict inequality $|P(\x)|>\dots$ can be
replaced by a non-strict one $|P(\x)|\ge\dots$ because
the volume of any level set of a non-constant polynomial is $0$.
$\square$

\medskip
Let now $\Cal F$ be a convex set in $\R^n$ of volume $\Vol (\Cal F)=1$,
and let $P$ be any (non-constant) polynomial in $\R^n$ of degree $d$.
Let $M(P)$ be the unique positive number such that
$$
\vol\{\x\in \Cal F\,:\, |P(\x)|\ge M(P)\}=1/e.
$$
\tit{Distribution inequalities:}
{\it For every} $\la>1$,
$$
\vol\{\x\in \Cal F\,:\,|P(\x)|>(A\la)^d M(P)\}\le e^{-\la}
$$
{\it and}
$$
\vol\{\x\in \Cal F\,:\,|P(\x)|<(A\la)^{-d} M(P)\}\le \frac1\la.
$$
\tit{Proof:}
The first inequality is just the comparison lemma applied to $c=M(P)$.
To obtain the second one, let us denote the volume on the left by $V$.
According to the comparison lemma applied to $c=(A\la)^{-d}M(P)$, we
have
$$
1/e\le (1-V)^{\la}
$$
and thereby
$$
V\le 1-e^{-1/\la}\le \frac 1\la.
\qquad\qquad\qquad \square
$$

\tt{Remarks:}
The first distribution inequality (basically due to Bourgain [DI1]) can be
viewed as (a kind of) {\it concentration phenomenon}. 
The second distribution inequality resembles a lot
the classical Remez estimate (R): the only difference is that instead of
the maximum over the {\it entire} set $\Cal F$, we have the {\it ``median"}
$M(P)$ on the left hand side. We want to emphasize here that
the comparison lemma and both distribution inequalities
are derived directly from the one-dimensional Remez inequality and, thereby,
remain valid (together with all their corollaries below)
for an arbitrary function (or class of functions) for which
the one-dimensional result holds. For instance, instead of polynomials of
degree $d$, we may consider exponential polynomials of order $d$, i.e.,
functions of the kind
$$
P(\x)=\sum_{k=1}^d c_k e^{i\langle\x_k,\x\rangle}
$$
with $c_k\in\C$, $\x_k\in\R^n$, for which the Remez inequality
(known in this case as Turan's lemma) holds with
$A=316$, say.

It may also be worth mentioning that replacing the
somewhat loose inequality (R) by the sharp one-dimensional
Remez estimate coming from the
consideration of Chebyshev polynomials, one can obtain the {\it sharp}
dimension-free comparison lemma
$$
\vol\{\x\in \Cal F\,:\, |P(\x)|\ge T_d(2\la-1)\, c \}\le
\bigl[\vol\{\x\in \Cal F\,:\,|P(\x)|\ge c\}\bigr]^{\la}
$$
and the corresponding distribution
inequalities
$$
\vol\{\x\in \Cal F\,:\,|P(\x)|>T_d(2\la-1)\, M(P)\}\le e^{-\la}
$$
and
$$
\vol\bigl\{\x\in \Cal F\,:\,|P(\x)|<\tfrac1{T_d(2\la-1)}\, M(P)
\bigr\}\le 1-e^{-\frac{1}{\la}}\,,
$$
where
$$
T_d(x)=\frac{1}{2}\left[
\bigl(x+\sqrt{x^2-1}\bigr)^d+\bigl(x-\sqrt{x^2-1}\bigr)^d
\right]
$$
is the Chebyshev polynomial of degree $d$.

\tt{Digression: estimates for average values via distribution functions
}
Since in what follows
we shall have to calculate a few integrals and averages
of real valued functions
using the estimates for their distribution functions,
let us remind the reader the corresponding general formulae.

Let $\Cal X$ be a measure space with measure $\mu$. Let $g:\Cal X\to \R$.
Let $\Cal Y$ be a measurable subset of $\Cal X$.
We want to construct a formula that would allow us to evaluate
the integral $\int_{\Cal Y}g\,d\mu$ or the average value
$\dsize\langle g\rangle\ci {\Cal Y}:=
\frac{1}{\mu(\Cal Y)}\int_{\Cal Y}g\,d\mu$
of the function
$g$ over the set $\Cal Y$
using only the information about measures of sets
of the kind $\{\x\in\Cal X\,:\, g(\x)>t\}$.

Fix some
``floor level" $L\in \R$ and consider the function $g\ui +:=\max(g,L)$.
We have
$$
g\ui+(\x)=L+\int_{L}^{g\ui+(\x)}dt.
$$
Hence,
$$
\int_{\Cal Y} g\,d\mu\le
\int_{\Cal Y} g\ui+\,d\mu\le
L\mu(\Cal Y)+
\int_L^{+\infty} \mu\{\x\in \Cal X\,:\, g(\x)>t\}\,dt\,,
$$
and, finally,
$$
\langle g\rangle\ci {\Cal Y}\le
L+\frac{1}{\mu(\Cal Y)}
\int_L^{+\infty} \mu\{\x\in \Cal X\,:\, g(\x)>t\}\,dt.
$$

For practical computations, we shall need the following modifications
of these estimates. Let $\f$ be any smooth, increasing to $+\infty$
function on $(0,+\infty)$.
Let $\La$ belong to
the domain of $\f$ and let $L=\f(\La)$. Making the change of variable
$t=\f(\la)$, we can rewrite
the above estimates as
$$
\aligned
\int_{\Cal Y} g\,d\mu
&\le
\f(\La)\mu(\Cal Y)+
\int_\La^{+\infty} \mu\{\x\in \Cal X\,:\, g(\x)>\f(\la)\}\f'(\la)\,dt
\\
\langle g\rangle\ci {\Cal Y}
&\le
\f(\La)+
\frac{1}{\mu(\Cal Y)}
\int_\La^{+\infty} \mu\{\x\in \Cal X\,:\, g(\x)>\f(\la)\}\f'(\la)\,d\la.
\endaligned
\tag $*$
$$

\tit{Estimates for $L^q$-norms:}
We shall need the following trivial observation:
for every $\sigma\ge 1$,
$$
\multline
1+\sigma\int_0^\infty \la^{\sigma-1}e^{-\la}\,d\la=
1+2^{\sigma-1}\sigma^{\sigma}\int_0^\infty
\left[\frac{\la}{2\sigma}\right]^{\sigma-1}e^{-\la}\,d\la
\\
\le 1+2^{\sigma-1}\sigma^\sigma\int_0^\infty e^{-\la/2}d\la
=1+(2\sigma)^\sigma\le (3\sigma)^\sigma.
\endmultline
$$
Let now $q\ge\frac{1}{d}$.
Applying the estimates ($*$) with $\Cal X=\Cal Y=\Cal F$,
$\mu=\vol$, $g=\Bigl[\dfrac{|P|}{A^d M(P)}\Bigr]^q$, $\f(\la)=\la^{qd}$,
$\La=1$ and
using the estimate $\mu\{g>\f(\la)\}\le e^{-\la}$
(which is equivalent to the first distribution inequality),
we get
$$
\int_{\Cal F} \left[\frac{|P|}{A^d M(P)}\right]^q\le
1+qd\int_1^\infty \la^{qd-1}e^{-\la}d\la\le (3qd)^{qd}.
$$
Therefore,
$$
\|P\|\ci{L^q(\Cal F)}\le (3Aqd)^d M(P)\quad \text{ for every }q\ge
\frac{1}{d}\,.
$$
Using the monotonicity of the function $q\to \|P\|\ci{L^q(\Cal F)}$, we
immediately derive from here that
$$
\|P\|\ci{L^q(\Cal F)}\le (3A)^d M(P) \quad \text{ for every }q\le
\frac{1}{d}\,.
$$
\tit{Estimates for $L^{-q}$-norms:}
Let $0<q<\frac{1}{d}$.
Applying the estimates ($*$) with $\Cal X=\Cal Y=\Cal F$,
$\mu=\vol$, $g=\Bigl[\dfrac{|P|}{A^{-d}M(P)}\Bigr]^{-q}$,
$\f(\la)=\la^{qd}$,
$\La=1$ and
using the estimate $\mu\{g>\f(\la)\}\le \dfrac 1\la$
(which is equivalent to the second distribution inequality), we get
$$
\int_{\Cal F} \left[\frac{A^d |P|}{M(P)}\right]^{-q}\le
1+qd\int_1^\infty \la^{qd-1}\frac{1}{\la}d\la = \frac{1}{1-qd}\,.
$$
Therefore,
$$
\|P\|\ci{L^{-q}(\Cal F)}\ge A^{-d}(1-qd)^{1/q}M(P)
\quad \text{ for every }0< q < \frac{1}{d}\,.
$$
\tit{The geometric mean:}
The above inequalities immediately imply that
$$
(eA)^{-d}M(P)\le \|P\|\ci{L^0(\Cal F)}\le (3A)^d M(P).
$$
\tit{Inverse H\"older inequalities:}
We shall start with the following simple
\tit{Observation:}
Let $A_+$ be the best constant such that
$$
\vol\{|P|\ge (A_+\la)^dM(P)\}\le e^{-\la} \quad\text {for all }\la\ge 1.
$$
Let $A_-$ be the best constant such that
$$
\vol\{|P|< (A_-\la)^{-d} M(P)\}\le \frac1{\la} \quad\text {for all }\la\ge
1.
$$
Then $A_+A_-\le A$.
\tit{Proof of the observation:}
Let $0<a<A_-$. According to the definition of $A_-$, there exists
$\la\cci-\ge 1$
such that
$$
\vol\{|P|< (a\la\cci-)^{-d} M(P)\}\ge \frac1{\la\cci-}
$$
and, thereby,
$$
\vol\{|P|\ge (a\la\cci-)^{-d} M(P)\}\le 1-\frac1{\la\cci-}\le
e^{-\frac{1}{\la\cci-}}.
$$
Then, for every $\la\ge 1$, we have
$$
\multline
\vol\{|P|\ge (\tfrac Aa\la)^d M(P)\}=
\vol\{|P|\ge (A[\la\la\cci-])^d
(a\la\cci-)^{-d} M(P)\}
\\
\le
\vol\{|P|\ge (a\la\cci-)^{-d} M(P)\}^{\la\la\cci-}\le
\Bigl[e^{-\frac{1}{\la\cci-}}\Bigr]^{\la\la\cci-}=e^{-\la}
\endmultline
$$
according to the comparison lemma applied with $c=(a\la\cci-)^{-d}$.
Thus, $A_+\le\frac{A}{a}$ and, since $a<A_-$ was arbitrary, we are done.
$\square$

\smallskip
Applying the above estimates for the $L^q$ and $L^{-r}$-norms with
$A_\pm$ in place of $A$, we conclude that
$$
\|P\|\ci{L^q(\Cal F)}\cdot
\|1/P\|\ci{L^r(\Cal F)}\le
\frac{ (3A\max\{1,qd\})^d }{(1-rd)^{1/r}}
\quad\text{ for all $q\ge 0$, $0\le r<1/d$.}
$$
\tit{The BMO - norm of $\log|P|$:}
We shall use the following definition of the $\BMO$-norm of a function
$u:\R^n\to\R$:
$$
\|u\|\ci{\BMO}=\sup\Sb \Cal F\subset \R^n \\
                      \Cal F\text{ is convex}\endSb\inf_{C\in\R}
\frac{1}{\vol(\Cal F)}\int_{\Cal F}|u-C|.
$$
Since the class of polynomials is dilation-invariant, it is
enough to obtain an estimate for convex sets $\Cal F$ of volume $1$.
Choosing $C=\log M(P)+d\,\frac{\log A_+-\log A_-}{2}$,
applying the estimates ($*$) with $\Cal X=\Cal Y=\Cal F$,
$\mu=\vol$, $g=\log|P|-C$, $\f(\la)=d\bigl[\log\la+\log (\sqrt A)\bigr]$,
$\La=1$, and
using the estimate $\mu\{g>\f(\la)\}\le
e^{-\la}+\dfrac{1}{\la}$
for the distribution function
(which is the combination of both estimates in the observation),
we get
$$
\int_{\Cal F}|\log |P|-C|\le
d\Bigl[\frac{\log A}{2}+\int_{1}^\infty\frac{1}{\la}\Bigl(e^{-\la}
+\frac{1}{\la}\Bigr)d\la\Bigr]\le
\frac{4+\log A}{2}\,d.
$$

\head
$\S 3$. Estimates for distribution of zeroes of ``random" analytic
functions
\endhead

\tt{An estimate for the averages of $\log|P|$ over subsets of a
compact convex set:}
\noindent The purpose of this subsection is to prove the following

\tt{Claim:}
{\it Let $\Cal F\subset \R^n$ be a compact convex set and
let $P:\R^n\to \R$ be a polynomial of degree $d$. Then for any
measurable $\Cal E\subset \Cal F$,}
$$
\bigl|
\langle \log|P|\rangle\ci {\Cal E} -
\langle \log|P|\rangle\ci {\Cal F}
\bigr|\le d \log\frac{e^2 A\vol(\Cal F)}{\vol(\Cal E)}\,,
$$
{\it where the averages are taken with respect to the $n$-dimensional
Lebesgue measure (volume) in $\R^n$ and $A$ is the constant in the
(one-dimensional) Remez inequality.}

\tt{Proof of the claim:}
Without loss of generality, we may assume that $\vol(\Cal F)=1$ and
$M(P)=1$.
Let, as before, $A_+$ and $A_-$ be the best constants in the inequalities
$$
\vol\{|P|\ge (A_+\la)^d\}\le e^{-\la}\le\frac{1}{\la} \quad(\la\ge 1)
$$
and
$$
\vol\{|P|< (A_-\la)^{-d}\}\le \frac1{\la} \quad(\la\ge 1).
$$
Taking
$\Cal X=\Cal F$, $\Cal Y=\Cal E$, $g=\log|P|$,
$\mu=\vol$, $\f(\la)=d(\log A_+ +\log\la)$, and using the inequality
$\mu\{g>\f(\la)\}\le \frac{1}{\la}$,
we
conclude that for every $\La\ge 1$,
$$
\langle\log|P|\rangle\ci {\Cal E} \le
d\Bigl[ \log A_+ +\log\La + \frac{1}{\vol(\Cal E)\La}\Bigr].
$$
Substituting $\La=\dfrac{1}{\vol(\Cal E)}$,  we get
$$
\langle\log|P|\rangle\ci {\Cal E} \le d\log\frac{e A_+}{\vol(\Cal E)}.
$$
Analogously, taking
$\Cal X=\Cal Y=\Cal F$, $g=-\log|P|$,
$\mu=\vol$, $\f(\la)=d(\log A_- +\log\la)$, we
conclude that for every $\La\ge 1$,
$$
\langle\log|P|\rangle\ci {\Cal F} \ge
-d\Bigl[ \log A_- +\log\La + \frac{1}{\La}\Bigr].
$$
Substituting $\La=1$,  we get
$$
\langle\log|P|\rangle\ci {\Cal F} \ge - d\log (e A_-).
$$
Combining these two estimates, we obtain
$$
\langle \log|P|\rangle\ci {\Cal E} -
\langle \log|P|\rangle\ci {\Cal F}
\le d\log\frac{e^2 A_+ A_-}{\vol(\Cal E)}
\le d\log\frac{e^2 A}{\vol(\Cal E)}.
$$
The inequality
$\dsize
\langle \log|P|\rangle\ci {\Cal E} -
\langle \log|P|\rangle\ci {\Cal F}
\ge- d\log\frac{e^2 A}{\vol(\Cal E)}
$
can be proved in a similar way.
$\square$

\tt{The Offord estimate:}
%
%
%
%
Fix some open domain $G\subset \C$ and consider a family of analytic
functions
$f(\x;\cdot):G\to\C$, where
$\x$ runs over some parameter set $\Cal X$ endowed with
a finite measure $\mu$. Let
$$
\nu_\x := \sum_{w:\, f(\x;w)=0} \delta_w
$$
be the {\it counting measure} of zeroes of the function $f(\x;\cdot)$
where $\delta_w$ stands for the Dirac measure at $w\in G$ and each
zero is counted with its multiplicity. For each $\x\in\Cal X$, the measure
$\nu_\x$ is a locally finite measure in $G$.

Consider the average measure
$$
\nu(U):=\frac{1}{\mu(\Cal X)}\int_{\Cal X}\nu_\x(U)\,d\mu(\x),\quad U\subset
G.
$$
The measure $\nu$ gives a ``typical" (average) distribution of
zeroes of the ``random" function $f(\x;\cdot)$ in $G$.
Let $\psi\in C_0^\infty(G)$ and let
$\la>0$.
Define the exceptional set $\Cal E_+=\Cal E_+(\psi,\la)$
by
$$
\Cal E_+(\psi,\la):=\Bigl\{\x\in\Cal X\,:\,
\int_G \psi\,d\nu_\x - \int_G \psi\,d\nu
\ge \la
\Bigr\}\,.
$$
Note that, since for each $\x\in\Cal X$,
the measure $\nu_\x$ is $\dfrac{1}{2\pi}$ times the
distributional Laplacian of the function $\log|f(\x;\cdot)|$,
we have
$$
\int_{G} \psi\,d\nu_\x =\frac{1}{2\pi}\int_{G}\Delta
\psi(z)\log|f(\x;z)|\,dm\cci 2(z)\,,
$$
where $m\cci 2$ is the area measure on the complex plane $\C$.
Averaging over $\Cal X$, we get
$$
\int_{G} \psi\,d\nu =\frac{1}{2\pi}\int_{G}\Delta
\psi(z)\langle\log|f(\cdot;z)|\rangle\ci{\Cal X}\,dm\cci 2(z).
$$
Averaging the difference of these identities with respect to the parameter
$\x$
over the set $\Cal E_+=\Cal E_+(\psi,\la)$, we obtain the inequality
$$
\multline
\la\le \frac{1}{2\pi}\int_G \Delta \psi(z)\cdot\bigl[
\langle \log|f(\cdot;z)|\rangle\ci{\Cal E_+}-
\langle \log|f(\cdot;z)|\rangle\ci{\Cal X}
\bigr]\,dm\cci 2(z)
\\
\le
\frac{1}{2\pi}\|\Delta\psi\|\ci{L^1(G)}\cdot \sup_{z\in G}
\left|
\langle \log|f(\cdot;z)|\rangle\ci{\Cal E_+}-
\langle \log|f(\cdot;z)|\rangle\ci{\Cal X}
\right|
\endmultline
$$
Almost exactly the same argument shows that the same inequality holds for
the
set
$$
\Cal E_-=\Cal E_-(\psi,\la):=\Bigl\{\x\in\Cal X\,:\,
\int_G \psi\,d\nu_\x - \int_G \psi\,d\nu
\le -\la
\Bigr\}\,.
$$
Combining these estimates with the claim, we obtain the following

\tt{Theorem (Offord's estimate):}
{\it If $\Cal X=\Cal F$ is a convex set in $\R^n$, $\mu$ is the Lebesgue
measure
in $\R^n$, and $f(\x;z)$ depends on $\x$ as a polynomial of degree $d$
for each $z$, then
$$
\frac{\Vol (\Cal E(\psi,\la))}{\Vol (\Cal F )}
\le
2Ae^2 \exp\biggl\{-\frac{2\pi \la}{d ||\Delta \psi||\ci {L^1(G)}}
\biggr\}\,,
$$
where}
$$
\Cal E(\psi,\la):=
\Cal E_+(\psi,\la)\cup
\Cal E_-(\psi,\la)=
\Bigl\{\x\in\Cal X\,:\,
\Bigl|\int_G \psi\,d\nu_\x - \int_G \psi\,d\nu\Bigr|
\ge \la
\Bigr\}\,.
$$

\tt{Corollary:}
Denote by $\Bbb D_r$ the disk of radius $r$ centered at the origin.
Let $G=\Bbb D_1$. We shall call a value
$\x \in \Cal F$ {\it exceptional} if the function
$f(\x;\cdot)$ does not vanish in $G$.
Let $\Cal E^*\subset \Cal F$ be the set of all exceptional values.
If $\vol(\Cal E^*)>0$, we can estimate the growth of the (average)
counting function $r \mapsto \nu
(\Bbb D_r)$\,:
$$
\nu (\Bbb D_r)
\le \frac{4 d}{1-r}
\log\frac{Ae^2\Vol (\Cal F)}{\Vol (\Cal E^*)}
\,, \qquad 0<r<1\,.
$$

\tt{Proof of the Corollary:}
Fix $r$ and choose a
test function $\psi (z) = \Psi (|z|)$ where $\Psi\in C_0^\infty [0,1)$,
$\Psi\ge 0$, and $\Psi(t)=1$ for $0\le t\le r$.
Let $\la:=\int_G\psi\,d\nu\ge \nu(\Bbb D_r)$. Note that for such choice of
$\la$, we obviously have $\Cal E^*\subset \Cal E_-(\psi;\la)$ and,
therefore,
$$
\frac{\Vol (\Cal E^*)}{\Vol (\Cal F)}
\le
Ae^2 \exp\biggl\{-\frac{2\pi \la}{d ||\Delta \psi||\ci{L^1(G)}} \biggr\}.
$$
We can rewrite it as
$$
\nu(\Bbb D_r)\le \la\le \frac{d}{2\pi}\|\Delta \psi\|\ci{L^1(G)}
\log\frac{Ae^2\Vol (\Cal F)}{\Vol (\Cal E^*)}\,.
$$
Note that
$$
\frac{1}{2\pi}\|\Delta \psi\|\ci{L^1(G)}
=\int_{r}^1 |t\Psi''(t)+\Psi'(t)|\,dt.
$$
Choosing $\Psi$ sufficiently close to the quadratic spline whose second
derivative is
$-\frac{4}{(1-r)^2}$ between $r$ and $\frac{1+r}{2}$ and
$+\frac{4}{(1-r)^2}$ between $\frac{1+r}{2}$ and $1$, we observe that the
right hand side can always be made less than
$\frac{4}{1-r}$. $\square$

\head
Appendix: Proof of the lemma
\endhead

Before starting the proof, we will make several simple observations about
numerical inequalities that we shall use in the course of the proof.

\tit{Observation 1:}
{\it For all} $X>0$, $Y\ge 0$,
$$
(X+Y)^\la \ge X\bigl(X+\tfrac{\la}{\la-1}Y\bigr)^{\la-1}.
$$

Indeed, we have an identity for $Y=0$, and, obviously, for each $Y\ge 0$,
$$
\frac{\dd}{\dd Y}\log(\LHS)=\frac{\la}{X+Y}\ge
\frac{\la}{X+\frac{\la}{\la-1}Y}=\frac{\dd}{\dd Y}\log(\RHS)
$$
where, as usual, $\operatorname{L(R)HS}$ stands for the Left (Right) Hand
Side of the inequality. $\square$
\tit{Observation 2:}
{\it If the inequality $(X+Y)^\la\ge X(X+Z)^{\la-1}$ holds for some
$X>0$, $Y,Z\ge 0$, then for each $T\ge 0$,}
$$
(X+Y+T)^\la\ge X\bigl(X+Z+\tfrac{\la}{\la-1}T\bigr)^{\la-1}.
$$

Indeed, if $Z\ge Y$, we may repeat the proof of
Observation 1 with $\frac{\dd}{\dd T}$ instead of $\frac{\dd}{\dd Y}$.
If $Z<Y<\frac{\la}{\la-1}Y$, then the desired inequality immediately
follows from Observation~1. $\square$

\tit{Observation 3:}
{\it If $(X+Y)^\la \ge X(X+Z)^{\la-1}$ for some $X,Y,Z>0$,
then
$$
(x+Y)^\la \ge x(x+Z)^{\la-1}\quad\text{ for all $x\in [0,X]$}.
$$
}

This is the least trivial of our observations.
Rewrite the inequality in the form
$$
\frac{x}{x+Y}\le \left[\frac{x+Y}{x+Z}\right]^{\la-1}
$$
which is equivalent to
$$
\left[\frac{x}{x+Y}\right]^{-\frac{1}{\la-1}}
\ge \frac{x+Z}{x+Y}
$$
Denote $\beta:=\frac{1}{\la-1}$,
$\theta:=\frac{x}{x+Y}$, $\Theta:=\frac{X}{X+Y}$.
Then $\frac{x+Z}{x+Y}=\frac{Z}{Y}-(\frac{Z}{Y}-1)\theta=L(\theta)$ is a
linear function.
We want to show that if the inequality $\theta^{-\beta}\ge L(\theta)$ holds
at $\theta=\Theta$, then it holds on the entire interval $[0,\Theta]$.
The desired inequality obviously holds for $\theta$ sufficiently
close to $0$. Therefore, if it were false for at least one
$\theta\in[0,\Theta]$, the graphs of functions $\theta^{-\beta}$ and
$L(\theta)$ would intersect at at least two points on the interval
$(0,\Theta]$. Since they also intersect at $\theta=1$, we would then
have at least three points common for a convex curve (the graph of
$\theta^{-\beta}$) and a line (the graph of $L(\theta)$), which is
impossible. $\square$

\tit{Observation 4:}
{\it If $(X+Y)^\la \ge X(X+Y+Z)^{\la-1}$ for some
$X,Y,Z>0$, then}
$$
(x+y)^\la \ge x(x+y+z)^{\la-1}
\quad
\text{ for all
$x\le X$, $y\ge Y$, $z\le Z$}.
$$

Indeed, we obviously can replace $Z$ by $z$. After that, Observation 3
allows
us to change $X$ to $x$. It remains to observe that for fixed $x$ and $z$,
$$
\frac{\dd}{\dd y}\log(\LHS)=\frac{\la}{x+y}\ge\frac{\la-1}{x+y+z}=
\frac{\dd}{\dd y}\log(\RHS). \qquad \qquad \square
$$

\medskip
Now we are ready to start proving the lemma.
Since the problem is invariant
with respect to linear change of variable, we may assume that $I=[0,1]$
We may also assume without loss of generality that the function $f$ is
continuous, {\it strictly} logarithmically concave and satisfies
$f(0)=f(1)=0$ (if it isn't
so, just consider the family of functions
$f_\e(x)=\bigl[x(1-x)\bigr]^\e f(x)$, apply the statement to each of them,
and pass to the limit as $\e\to 0$).

Clearly, $E\ci{\la, I}$ is a closed set. If $E\ci{\la,I}$ is empty,
there is nothing to prove.
Otherwise, $(0,1)\setminus E\ci{\la,I}=\cup_j I_j$ where $I_j$ are disjoint
open
intervals each of which is shorter than the entire interval $(0,1)$.
Consider one of these intervals $I_j=(a,b)$. We shall call it
{\it regular} if either $a>0$ and $f$ is decreasing on $(a,b)$, or $b<1$
and
$f$ is increasing on $(a,b)$. Otherwise we shall call the interval
$I_j$ {\it exceptional}.
Clearly, there may be not more than one exceptional interval. If such an
interval exists, we shall assign the index $0$ to it. Let $E_j=E\cap I_j$.
We claim that for each regular interval, one has
$$
\int_{E_j}f\ge \frac{\la-1}\la\int_{I_j}f.
$$
Indeed, if, say, $I_j=(a,b)$ and $a>0$, then $a\in E\ci{\la,I}$ and,
thereby, $|E_j\cap(a,t)|\ge \frac{\la-1}\la
|(a,t)|$ for each $a<t<b$, which, together with the fact that
$f$ is decreasing on $(a,b)$, is enough to ensure the desired estimate.

If the exceptional interval is absent, the inequality of the
lemma is quite easy to prove. Indeed, it is equivalent to the estimate
$$
\Bigl[\int_{E} f\Bigr]^\la\ge
\Bigl[\int_{E\cci{\la,I}} f\Bigr]\cdot
\Bigl[\int_{(0,1)} f\Bigr]^{\la-1}\,;
$$
i.e., to the inequality
$$
\Bigl[\int_{E\cci{\la,I}} f + \int_{\cup E_j}f\Bigr]^\la\ge
\Bigl[\int_{E\cci{\la,I}} f\Bigr]\cdot
\Bigl[\int_{E\cci{\la,I}} f + \int_{\cup I_j}f\Bigr]^{\la-1}.
$$
But $\int_{\cup I_j}f\le \frac{\la}{\la-1}\int_{\cup E_j}f$
and thereby the desired estimate follows from Observation 1.

Suppose now that $I_0=(a,b)$ is exceptional.
Without loss of generality we may
assume that $f(b)\le f(a)$ (otherwise
just make the change of variable $t\to 1-t$, which leaves the problem
invariant). Note that this automatically implies that $a>0$
because otherwise we would have $f(b)\le f(a)=f(0)=0$,
which, since the function $f$ is strictly positive
on $(0,1)$, would imply that $b=1$, $I\ci 0=(0,1)$, and, finally,
that $E\ci{\la,I}$ is empty.

If $f(b)<f(a)$, let $c\in(a,b)$ be the
(unique) point such that $f(a)=f(c)$.
We are going to slightly modify the portion $E_0$ of the set $E$.
Observe again that, since $a\in E\ci{\la,I}$, $|E\cap(a,c)|\ge
\frac{\la-1}\la |(a,c)|$.
Take an arbitrary portion of $E\cap(a,c)$ of measure
$|E\cap(a,c)|-\frac{\la-1}\la|(a,c)|$
and replace it by a set of equal measure on
$(c,b)$ using the points of $(c,b)\setminus E$ as close to the left end $c$
as possible. If the measure of the entire set $(c,b)\setminus E$ is too
small, just fill the entire interval $(c,b)$ and forget about lost measure.
Let $E'$ be the resulting set.
We claim that
$$
|E'\cap(c,t)|\ge \frac{\la-1}\la|(c,t)|
$$
for all $c<t<b$.
Indeed, the portion $E'\cap(c,b)$ of the modified set $E'$ starts with an
interval. As long as $t$ stays within this interval, there is nothing to
prove. As soon as $t$ leaves this interval, the length of the intersection
$E'\cap(a,t)$ coincides with the length of the intersection $E\cap(a,t)$
and therefore is not less than $\frac{\la-1}\la |(a,t)|$.
But we also have $|E'\cap(a,c)|=\frac{\la-1}\la|(a,c)|$,
so we should have the desired inequality for the remaining portion.
Also, $f$ obviously decreases on $(c,b)$. So, we may treat the interval
$(c,b)$ as a regular interval and to restrict our attention to $(a,c)$.

If we originally had the identity $f(a)=f(b)$, this construction reduces
to denoting the point $b$ by the letter $c$ and replacing the part
$E\cap (a,b)$ of the set $E$ by its arbitrary subset of measure
$\frac{\la-1}\la|(a,b)|$.

On $(a,c)$, let us modify the set $E'$ even further. Namely, let us replace
the
corresponding portion of $E'$ by the level set of $f$ of measure
$\frac{\la-1}\la|(a,c)|$ containing the small values of the function.
Clearly, such modifications only decrease the integral of the function $f$
over the set that undergoes them, so we have $\int_{E'}f\le\int_E f$.
Thus, it will suffice to prove the inequality of the lemma
with $\int_E f$ replaced by $\int_{E'}f$.

Now let us look at the picture we have obtained. We have one exceptional
interval $I'_0=(a,c)$ such that $f(a)=f(c)$,
$|E'_0|=\frac{\la-1}\la|(a,c)|$,
and $E'_0$
is a level set of $f$ on $I'_0$ containing the small values of the function.
We have also some regular intervals $I'_j$ (original regular intervals plus,
maybe,
$(c,b)$\ ) satisfying $\int_{E'_j}f\ge\frac{\la-1}\la\int_{I'_j}f$ for each
$j$.
We need to prove the estimate
$$
\Bigl[\int_{E\cci{\la,I}} f + \int_{E'_0}f+ \int_{\cup_{j>0}
E'_j}f\Bigr]^\la\ge
\Bigl[\int_{E\cci{\la,I}} f\Bigr]\cdot
\Bigl[\int_{E\cci{\la,I}} f +\int_{I'_0}f+\int_{\cup_{j>0} I'_j
}f\Bigr]^{\la-1}.
$$
Using Observation 2, we see that it is enough to prove that
$$
\Bigl[\int_{E\cci{\la,I}} f + \int_{E'_0}f\Bigr]^\la\ge
\Bigl[\int_{E\cci{\la,I}} f\Bigr]\cdot
\Bigl[\int_{E\cci{\la,I}} f +\int_{I'_0}f\Bigr]^{\la-1}.
$$
Observation 3 allows us to extend the set $E\cci{\la,I}$
in the last inequality to the entire
set $(0,1)\setminus I'_0$. Let now $|I'_0|=\la m$ and let $f^*$ be the
decreasing
rearrangement of $f$ on $(0,1)$. It is obviously decreasing and
logarithmically concave. We need to prove the inequality
$$
\Bigl[\int_{\la m}^1 f^*+\int_{m}^{\la m}f^*\Bigr]^\la\ge
\Bigl[\int_{\la m}^1 f^*\Bigr]
\Bigl[\int_{\la m}^1 f^*+\int_m^{\la m}f^* + \int_0^m f^*\Bigr]^{\la-1}.
$$
According to Observation 4, if we modify $f^*$ in such a way
that simultaneously the integrals $\int_0^m f^*$ and $\int_{\la m}^1 f^*$
become
bigger while the integral $\int_m^{\la m} f^*$ becomes smaller, we shall get
a
harder inequality to prove. Such modification can be done by replacing
$\log f^*$ by a linear function interpolating it at the points $m$ and $\la
m$.
Using Observation 3 once more, we see that we may extend the
integration to the entire
right semi-axis. Finally, we need to prove that if $f^*$ is a decreasing
exponential function, then
$$
\Bigl[\int_m^\infty f^*\Bigr]^\la\ge
\Bigl[\int_{\la m}^\infty f^*\Bigr]\cdot
\Bigl[\int_0^\infty f^*\Bigr]^{\la-1}.
$$
But this is an identity!
$\square$

\tit{Acknowledgements:} The authors thank Efim Gluskin, Vitali Milman and
Leonid Polterovich for useful discussions.

\define\ba#1#2#3#4#5#6#7
{\par\smallskip\noindent{{\bf #1.}\,\,{\rm #2,}\,\,{\it #3}\,\,\,\,}
\par\noindent{#4}\,\,{\bf #5}\,(#6)\,#7.}

\bigskip
\centerline{\bf Related Literature}

\bigskip

\noindent{\bf Brunn-Minkowski theorem}

\smallskip

\ba {BM1}{Yu. D. Burago and V. A. Zalgaller}{Geometric Ineqalities}
{Springer-Verlag}{}{1988}

\bigskip

\noindent{\bf Needle decomposition}

\smallskip

\ba {ND1}{M. Gromov and V. Milman}{Generalization of the spherical
isoperimetric inequality
to uniformly convex Banach spaces}{Compositio Math.}{62}{1987}{263-282}

\ba {ND2}{L. Lov\'asz and  M. Simonovits}{Random walks in a convex body
and an
improved volume algorithm}{Random Structures and
Algorithms}{4}{1993}{359--412}

\ba {ND3}{R. Kannan, L. Lov\'asz and M. Simonovits}
{Isoperimetric problem for convex bodies and a localization lemma}
{Discrete and Comput. Geometry}{13}{1995}{541--559}

\bigskip

\noindent{\bf Remez inequality}

\smallskip

\ba {RI1} {E. J. Remez}
{Sur une propri\'ete des polyn\^omes de Tschebycheff}
{Commun. Inst. Sci. Kharkov}{13}{1936}{93--95}

\ba {RI2}{R. M. Dudley and B. Randol}{Implications of pointwise bounds on
polynomials}{Duke Math. J.} {29}{1962}{455--458}

\ba {RI3}{Yu. Brudnyi and M. Ganzburg}{One extremal problem for
polynomials
of $n$ variables}{Izv.Akad.Nauk SSSR (Mat)} {37}{1973}{344--355, (in
Russian)}

\bigskip

\noindent{\bf Tur\'an lemma}

\smallskip

\ba {TL1}{F. Nazarov}
{Local estimates for exponential polynomials and their
applications to the uncertainty principle type results}
{Algebra i Analiz}{5}{1993}{3--66, (in Russian)}
\par\noindent English transl. St. Petersburg Math. J. (1994)

\bigskip

\noindent{\bf Dimension-free distribution inequalities}

\smallskip

\ba {DI1}{M. Gromov and V. Milman}{Brunn theorem and a concentration of
volume of convex bodies}
{GAFA Seminar Notes Tel Aviv University}{}{1983/84}{}

\ba {DI1}{J. Bourgain}{On the distribution of polynomials on high
dimensional
convex sets}{Lect. Notes in Math.}{1469}{1991}{127--137}

\ba {DI2}{A. Carbery and J. Wright}{Distributional and $L^q$ norm 
inequalities for polynomials over convex bodies in ${\Bbb R}^n$}
{Math. Res. Lett.}{8}{2001}{233--248}

\bigskip

\noindent{\bf Inverse H\"older inequality}

\smallskip

\ba {IH1}{D. Ullrich}
{Khinchin's inequality and the zeroes of Bloch functions}
{Duke Math. J.}{57}{1988}{519--535}

\ba {IH2}{V. Milman and A. Pajor}
{Cas limites dans des in\'egalit\'s du type de
Khinchine et applications g\'eom\'triques}
{C. R. Acad. Sci. Paris S\'er. I Math.}{ 308}{1989}{91--96}

\ba {IH3}{S. Favorov and E. Gorin}{Generalizations of Khinchin
inequality}{Theory Prob. Appl.} { 35} {1990}{766--771}

\ba {IH4}{S. Favorov}{A generalized Kahane-Khinchin inequality}{Studia
Mathematica} {130} {1998}{101-107}

\ba {IH5}{R. Lata\l a}
{On the equivalence between geometric and arithmetic means for
log-concave measures}{in: Convex geometric analysis (Berkeley, CA, 1996)
Math. Sci. Res. Inst. Publ., 34, Cambridge Univ. Press,
Cambridge}{}{1999}{123--127}

\ba {IH6}{O. Gu\'edon}{Kahane-Khinchine type inequalities for negative
exponents}{Mathematika} {46}{1999}{165--173}

\ba {IH7}{S. G. Bobkov}{Remarks on the growth of $L^p$ norms of
polynomials}
{in: "Geom. Aspects of Functional Analysis", Lecture Notes in
Math.}{1745}{2000}{27--35}

\bigskip

\noindent{\bf Offord's estimate}

\smallskip

\ba {OE1}{A. C. Offord}{The distribution of zeros of power series whose
coefficients
are independent random variables}{Indian J. Math.} {9} {1967}{175--196}

\ba {OE2}{M. Sodin}{Zeros of gaussian analytic functions}{Math. Res.
Lett.}{7}{2000}{371--381}

\enddocument